\documentclass[12pt]{article}

\usepackage{amsmath,amssymb,amsthm}
\usepackage[titletoc,title]{appendix}
\usepackage{float}
\usepackage{bm}
\usepackage{latexsym}
\usepackage{cases}
\usepackage[utf8]{inputenc} 
\usepackage[T1]{fontenc}

\newtheorem{dfn}{Def}[section]
\newtheorem{thm}[dfn]{Theorem}
\newtheorem{prop}[dfn]{Proposition}
\newtheorem{lem}[dfn]{Lemma}

\allowdisplaybreaks[3]

\makeatletter

\newcounter{myc}[section]

\newcounter{mycc}[section]

\def\id#1{\def\@id{#1}}
\def\department#1{\def\@department{#1}}
\def\superadvisor#1{\def\@superadvisor{#1}}

\makeatletter
    
    \@addtoreset{equation}{section}
\makeatother

\begin{document}
\title{\large
\bf{On $L^1$ estimates of solutions\\ of compressible viscoelastic system}}
\author{
\large{Yusuke Ishigaki}}
\date{\small
Department of Mathematics, \\
Tokyo Institute of Technology, \\
Tokyo 152-8551, Japan, \\
e-mail:ishigaki.y.aa@m.titech.ac.jp\\
}

\maketitle
\begin{abstract}
We consider the large time behavior of solutions of compressible viscoelastic system around a motionless state in a three-dimensional whole space. We show that if the initial preturbation belongs to $W^{2,1}$, and is sufficiently small in $H^4\cap L^1$, the solutions grow in time at the same rate as $t^{\frac{1}{2}}$ in $L^1$ due to diffusion wave phenomena of the system caused by interaction between sound wave, viscous diffusion and elastic wave.   
\end{abstract}

\section{Introduction}

This paper studies the initial value problem of compressible viscoelastic system   
\begin{gather}
\partial_t\rho+\mathrm{div} (\rho v)=0, \label{system-eq1} \\[1.0ex]
\rho(\partial_t v+ v\cdot\nabla v)
-\nu{\Delta}v -(\nu+\nu')\nabla \mathrm{div} v
+\nabla  P(\rho)  = \beta^2\mathrm{div}(\rho F{}^\top\! F), \label{system-eq2} \\[1.0ex]
\partial_t F + v\cdot\nabla F =(\nabla v)F, \label{system-eq3} \\[1.0ex]
(\rho, v, F)|_{t=0}=(\rho_0,v_0,F_0). \label{system-IC}
\end{gather}
in the whole space $\mathbb{R}^3$.
Here $\rho=\rho(x,t)$, 
$v ={}^\top (v^1(x,t), v^2(x,t), v^3(x,t))$, and $F=(F^{jk}(x,t))_{1\leq j,k\leq 3}$ 
 denote the unknown density, velocity field, and deformation tensor, respectively, at position $x\in \mathbb{R}^3$ and time $t\geq 0$;
$P=P(\rho) $ is the given pressure; $\nu$ and $\nu'$ are the viscosity coefficients satisfying
\[
\nu>0,~2\nu+3\nu'\geq0;
\]
$\beta>0$ is the propagation speed of elastic wave. 
We assume that $P'(1)>0$, and set $\gamma=\sqrt{P'(1)}$. 

We also impose the following conditions for $\rho_0$ and $F_0$
\begin{gather}
\rho_0\mathrm{det}F_0=1, \label{initialcond21}\\
\sum_{m=1}^3(F_0^{ml}\partial_{x_m}F_0^{jk}-F_0^{mk}\partial_{x_m}F_0^{jl})=0,
~j,k,l=1,2,3, \label{initialcond22}\\
\mathrm{div}(\rho_0 {}^\top\! F_0)=0.\label{initialcond23}
\end{gather}
It follows from [5, Appendix A] and [18, Proposition.1] that the quantities $(1.5)$ and $(1.6)$ are invariant for $t\geq0$:
\begin{gather}
\rho\mathrm{det}F=1, \label{constraint1.4.1}\\
\sum_{m=1}^3(F^{ml}\partial_{x_m}F^{jk}-F^{mk}\partial_{x_m}F^{jl})=0,~j,k,l=1,2,3.
\label{constraint1.4.2}
\end{gather}
Here the constraint $(1.8)$ means compressibility of fluid and the constraint $(1.9)$ called the Piola's formula is derived from a certain symmetric property of the first order derivatives of $F$ in the Lagrangian coordinates. 
Furthermore we see from [7. Appendix A] and [8. Appendix A] that the constraints $(1.8)$ and $(1.9)$ lead to the time invariance of the quantity $(1.7)$: 
\begin{align}
\mathrm{div}(\rho{}^\top\! F)=0. \label{constraint1.4.3}
\end{align}
\vspace{1ex}

The purpose of this paper is to deduce the estimate of the $L^1$ norm of solutions of the problem $(1.1)$--$(1.7)$ around a motionless state $(\rho,v,F)=(1,0,I)$. Here $I$ is the $3\times3$ identity matrix.
\vspace{1ex}

The system $(1.1)$--$(1.3)$ is derived from a motion of compressible viscoelastic fluid in the macroscopic scale under the Hookean li3near elasticity by the variational settings. We refer to [1, 13, 21] for more physical details.  We can classify the system $(1.1)$--$(1.3)$ in a quasilinear parabolic-hyperbolic system since the system $(1.1)$--$(1.3)$ consists of the compressible Navier-Stokes equations and a first order hyperbolic system for $F$.  
\vspace{1ex}

In the case $\beta=0$, the large time behavior of the solutions around $(\rho,v)=(1,0)$ has been investigated so far. 
In particular, if we set $\beta=0$ formally, the system $(1.1)$--$(1.3)$ becomes the usual compressible Navier-Stokes equation
\begin{gather}
\partial_t\rho+\mathrm{div} (\rho v)=0, \label{CNS-eq1} \\[1.0ex]
\rho(\partial_t v+ v\cdot\nabla v)
-\nu{\Delta}v -(\nu+\nu')\nabla \mathrm{div} v
+\nabla  P(\rho)  = 0, \label{CNS-eq2} \\[1.0ex]
(\rho, v)|_{t=0}=(\rho_0,v_0). \label{CNS-IC}
\end{gather}

Matsumura and Nishida [15] showed the global existence of the solutions of the problem $(1.11)$--$(1.13)$ provided that the initial perturbation is sufficiently small in $H^3\cap L^1$, and derived the decay estimate:
$$
\|\nabla^k(\phi(t),m(t))\|_{L^2}\leq C(1+t)^{-\frac{3}{4}-\frac{k}{2}}, ~k=0,1,
$$ 
where $(\phi,m)=(\rho-1,\rho v)$.
Hoff and Zumbrun [2] derived the following $L^p$ decay estimates and asymptotic properties in $\mathbb{R}^n,~n\geq 2$:
\begin{gather*}
\|(\phi(t),m(t))\|_{L^p}\leq
\begin{cases}
C(1+t)^{-\frac{n}{2}\left(1-\frac{1}{p}\right)-\frac{n-1}{4}\left(1-\frac{2}{p}\right)}L(t), & 1\leq p < 2, \\
C(1+t)^{-\frac{n}{2}\left(1-\frac{1}{p}\right)}, & 2\leq p \leq\infty,
\end{cases}\\[1ex]
\begin{array}{l}
\left\|\left((\phi(t),m(t))-\left(0,\mathcal{F}^{-1}\left(e^{-\nu|\xi|^2 t}\hat{\mathcal{P}}(\xi)\hat{m}_0\right)\right)\right)\right\|_{L^p} \\[1.5ex]
\hspace{1.5ex}
\leq C(1+t)^{-\frac{n}{2}\left(1-\frac{1}{p}\right)-\frac{n-1}{4}\left(1-\frac{2}{p}\right)}, 
\end{array}
~2\leq p \leq \infty,
\end{gather*}
provided that the initial perturbation is sufficiently small in $H^4\cap L^1$, where $L(t)=\log(1+t)$ when $n=2$, and $L(t)=1$ when $n\geq3$.  Here $\hat{\mathcal{P}}(\xi)=I-\frac{{\xi}^\top \xi}{|\xi|^2},~\xi\in\mathbb{R}^n$. According to [10], the solution of the linearized system is decomposed as the sum of two terms, one is the incompressible part given by $\mathcal{F}^{-1}\left(e^{-\nu|\xi|^2 t}\hat{\mathcal{P}}(\xi)\hat{m}_0\right)$ which behaves pure diffusively, and the other is the compressible part  $(\phi(t),m(t))-\left(0,\mathcal{F}^{-1}\left(e^{-\nu|\xi|^2 t}\hat{\mathcal{P}}(\xi)\hat{m}_0\right)\right)$ containing the diffusion wave which stands for the convolution of the heat kernel and the fundamental solution of the wave equation with sound speed $\gamma$. The authors of [2] revealed the hyperbolic aspect of the system $(1.11)$--$(1.12)$ by proving that the large time behavior of the compressible part is different from the heat kernel in $L^p$ as $t\to\infty$, except the case $p=2$. See also [11] for the linearized problem.   
\vspace{1ex}

In the case $\beta>0$, the mathematical analysis of solutions of the initial value problem $(1.1)$--$(1.7)$ around the motionless state have been developed so far. The local existence of the strong solution of the problem $(1.1)$--$(1.7)$ was guaranteed by Hu and Wang [4]. The global existence of the strong solution of the initial value problem $(1.1)$--$(1.7)$ was proved by Hu and Wang [5], Qian and Zhang [18], and Hu and Wu [6], provided that the initial perturbation $(\rho_0-1,v_0,F_0-I)$ is sufficiently small. Hu and Wu [6] showed that if the initial perturbation $(\rho_0-1,v_0,F_0-I)$ belongs to $L^1\cap H^3$, the $L^p$ decay estimates hold for the case $2\leq p\leq6$:
\begin{align}
\|u(t)\|_{L^p}\leq C(1+t)^{-\frac{3}{2}\left(1-\frac{1}{p}\right)}, \label{LpdecayestCh1}
\end{align}
by using the Fourier splitting method and the Hodge decomposition. 
Here $u(t)=(\phi, w, G)=(\rho-1,v,F-I)$. Moreover, the authors also derived the lower $L^2$ estimate
\begin{align}
\|u(t)\|_{L^2}\geq c(1+t)^{-\frac{3}{4}},~t\gg 1, \label{sharpdecayestCh1}
\end{align}
provided that the following conditions satisfy in the low frequency part $|\xi|\ll1$:
\begin{align}
|\hat{\phi}_0(\xi)|\geq c_0, |\hat{m}_0 (\xi)|+|\hat{\mathcal{G}}_0(\xi)-{}^\top\hat{\mathcal{G}}_0(\xi)| \ll |\xi|^{\eta_0}, \label{lowfreqHuWu}
\end{align}
where $(\phi_0,m_0,\mathcal{G}_0)=(\rho_0-1,\rho_0v_0,\rho_0F_0-I)$; $c$, $c_0$ and $\eta_0$ are positive constants independent of $\xi$ and $t$. Li,  Wei and Yao [12,~22] generalized the upper $L^p$ decay estimates $(1.14)$ to the case $2\leq p\leq \infty$, and obtained the $L^2$ decay estimates of higher order derivatives:  
\begin{align}
\|\nabla^k u(t)\|_{L^2}\leq C(1+t)^{-\frac{3}{4}-\frac{k}{2}},~k=0,1,\ldots,N-1,
\label{L2decayestCh1}
\end{align}
provided that $u_0=(\rho_0-1,v_0,F_0-I)$ belongs to $H^N,~N\geq3$, and is small in $L^1\cap H^{3}$. We also refer to [3, 14, 23] in recent progresses.

One problem which interests us is that the decay rates in $(1.14)$ reveal only the parabolic aspect of the system $(1.1)$--$(1.3)$; it would be desirable to establish decay estimates which reflect the hyperbolic aspect of the system $(1.1)$--$(1.3)$.
\vspace{1ex}

 In view of the results in [2], it is expected that the system $(1.1)$--$(1.3)$ has the diffusion wave phenomena affected by sound wave, viscous diffusion and elastic wave. For simplicity, let us consider the linearized system around $(1,0,I)$:
\begin{align}
\partial_t u+Lu=0. \label{LSintro}
\end{align}
Here $L$ denotes the linearized operator given by
\[
L=
\left(
\begin{array}{@{\ }cc@{\ }cc@{\ } }
0 & \mathrm{div}& 0\\
\gamma^2\nabla & -\nu\Delta-\tilde{\nu}\nabla\mathrm{div} & -\beta^2\mathrm{div} \\
0 & -\nabla & 0 
\end{array}
\right),
\]
where $\tilde\nu=\nu+\nu'$. We then see that the solenoidal part of the velocity $w_s=\mathcal{F}^{-1}(\hat{\mathcal{P}}(\xi)\hat{w})$ satisfies the following linear symmetric parabolic-hyperbolic system:
\[
\left\{
\begin{array}{l}
\partial_t w_s-\nu\Delta w_s-\beta\mathrm{div}\tilde{G}_s=0,\\
\partial_t \tilde{G}_s- \beta\nabla w_s=0,
\end{array}
\right.
\]
where $\tilde{G}_s=\beta\mathcal{F}^{-1}(\hat{\mathcal{P}}(\xi)\hat{G})$,
while the complimentary part $w_c=w-w_s$ solves the following strongly damped wave equation:  
\[
\partial_t^2 w_c-(\beta^2+\gamma^2)\Delta w_c-(\nu+\tilde\nu)\partial_t\Delta w_c=0.
\]
Owing to the results in [20], the large time behavior of the solution of $(1.18)$ becomes different to the case $\beta=0$ ([2,11]) since the additional hyperbolic aspect arises in the incompressible part due to elastic wave. As a result, the principal part of the linearized system $(1.18)$ can be identified as a system of the strongly damped wave equation. In [9], the hyperbolic aspect of the system $(1.1)$--$(1.3)$ is clarified by showing the following $L^p$ decay estimates:
\begin{align}
\|u(t)\|_{L^p}\leq C(1+t)^{-\frac{3}{2}\left(1-\frac{1}{p}\right)-\frac{1}{2}\left(1-\frac{2}{p}\right)}. \label{LpdecayestCh12}
\end{align}
for the case $1< p\leq\infty$, provided that the initial perturbation $(\rho_0-1,v_0,F_0-I)$ is small in $L^1\cap H^{3}$,
This improves the results in [6, 12]. 
\vspace{1ex}

The main difficulty of the mathematical analysis is nonlinearity of the constraints $(1.8)$--$(1.10)$. Therefore, straightforward application of the semigroup theory to the nonlinear problem does not seem valid. To bypass this difficulty, Hu and Wu [6] found that the behavior of $G$ is controlled by its skew-symmetric part $G-{}^\top G$ due to the constraints $(1.9)$ and $(1.10)$. This property leads to the global in time existence theorem. The authors of [6] next used the Helmholtz decomposition of $w$ and the skew-symmetric part of $G$ to obtain $(1.14)$ with $2\leq p\leq 6$, $(1.17)$ with $N=2$ and $(1.15)$.  In [9], the author relied on a material coordinate transform which makes the constraint $(1.10)$ a linear one to apply the analysis of the linearized problem to the nonlinear problem. Let us introduce a displacement vector $\tilde\psi=x-X$ as in [19, 21]: 
\[ 
\tilde\psi(x,t)=x-X(x,t).
\]
Here $x=x(X,t)$ is the material coordinate defined under the flow map
\[
\left\{
\begin{array}{l}
\displaystyle
\frac{\mathrm{d}x}{\mathrm{d}t}=v(x(X,t),t), \\[1ex]
x(X,0)=X,
\end{array}
\right.
\]
and $X=X(x,t)$ denotes the inverse of $x$.  In the continuum mechanics theory, $F$ is given by the Jacobi matrix of $x$ in the material coordinate. Then we notice that $F$ is written as
\begin{align}
F-I=\nabla\tilde{\psi}+h(\nabla\tilde{\psi}). \label{FIpsi}
\end{align} 
Here $h(\nabla\tilde{\psi})$ is a function satisfying  $h(\nabla\tilde{\psi})=O(|\nabla\tilde\psi|^2)$ for $|\nabla\tilde\psi|\ll1$.  The author next set the nonlinear transform 
\begin{align}
\psi=\tilde{\psi}-(-\Delta)^{-1}\mathrm{div} {}^\top (\phi\nabla\tilde{\psi}+(1+\phi) h(\nabla\tilde{\psi})). \label{NTrans}
\end{align}
 Here $(-\Delta)^{-1}=\mathcal{F}^{-1}|\xi|^{-2}\mathcal{F}$. By straight computation and the constraint $(1.10)$, we arrive at the linear condition $\phi+\mathrm{tr}(\nabla\psi)=\phi+\mathrm{div}\psi=0$ which makes the semigroup $e^{-tL}$ generated  by $-L$ tend to $0$ as $t\to\infty$ in $L^p,~p>\frac{5}{4}$. Furthermore, the decay estimate of the $L^p~(1< p\leq \infty)$ norm of $u=(\phi,w,G)$ is obtained from $\tilde{U}=(\phi,w,\nabla\psi)$. Consequently, the $L^p$ decay estimates of $u$ are obtained from the following integral equation
\[
\tilde{U}(t)=e^{-tL}\tilde{U}(0)+\int_0^t e^{-(t-s)L}N(\tilde{U}(s)) \mathrm{d}s,
\]
where $N(\tilde{U})=(N_1(\tilde{U}),N_2(\tilde{U}),N_3(\tilde{U}))$ is a nonlinearity satisfying $N_1(\tilde{U})+\mathrm{tr}N_3(\tilde{U})=0$. 
\vspace{1ex}

The another difficulty arises from the nonlinear transform $(1.21)$ containing the nonlocal operator $(-\Delta)^{-1}$. We note that the operator $\nabla(-\Delta)^{-1}\mathrm{div}=\mathcal{F}^{-1}\frac{{\xi}^\top \xi}{|\xi|^2}\mathcal{F}$ is not bounded from $L^p$ into $L^p$ if $p=1,\infty$ due to the Riesz operator. In case $p=\infty$, the above difficulty is avoided by using the Sobolev inequality and the Plancherel theorem, while the case $p=1$, it is expected that the solution of $L^1$ norm grows as $t\to\infty$ due to the diffusion phenomena, however, it remains open.
\vspace{1ex}

In this paper, we find a different approach of the reformulation to show the following $L^1$ estimate of $u$
\begin{align}
\|(\rho(t)-1, v(t), F(t)-I)\|_{L^1}\leq C(1+t)^{\frac{1}{2}},~t>0, \label{L1estCh1}
\end{align}   
provided that the initial perturbation $u_0=(\rho_0-1, v_0, F_0-I)$ belongs to $ W^{2,1}$, and is sufficiently small in $H^4\cap L^1$.
We also prove that if $(\rho_0,v_0,F_0)$ satisfies $(1.16)$, then the following lower $L^1$ estimate holds
\begin{align}
\|(\rho(t)-1, v(t), F(t)-I)\|_{L^1}\geq C(1+t)^{\frac{1}{2}},~t\gg 1, \label{L1sharpCh1}
\end{align}   
This indicates that the obtained rate $(1+t)^{\frac{1}{2}}$ in $(1.22)$ is sharp.
\vspace{1ex}

We give an outline of the proof of the main result. We first notice that the constraint $(1.8)$ is read as $\rho=\mathrm{det}F^{-1}$. Then, by employing $(1.20)$, we have
\begin{align}
\phi=-\mathrm{div}\tilde{\psi}+O(|\nabla\tilde{\psi}|^2),~\|\nabla\tilde{\psi}\|_{C(0,\infty;L^\infty)}\ll 1. \label{Ch1phipsiest}
\end{align}
Therefore, the behavior of $\phi$ can be handled by $-\mathrm{div}\tilde{\psi}$ under the small perturbation. For simplicity, we omit the tilde $\tilde{\cdot}$ of $\tilde{\psi}$ here.

We next consider the nonlinear problem for $U=(\tilde{\phi},w,\tilde{G})=(-\mathrm{div}\psi,w,\nabla\psi)$:
\begin{equation}
\left\{\label{system3Ch1}
\begin{array}{l}
\partial_t U+ LU= N(U),\\
\tilde{\phi} +\mathrm{tr}\tilde{G} =0,~\tilde{G}=\nabla\psi,\\
U|_{t=0}=U_0.
\end{array}
\right.
\end{equation}
where $N(U)=(N_1(U),N_2(U),N_3(U))$ is a nonlinearity such that $N_1(U)+\mathrm{tr}N_3(U)=0$. We see from $(1.20)$ and $(1.24)$ that the $L^1$ norm of $u=(\phi,w,G)$ is estimated by $U=(-\mathrm{div}\psi,w,\nabla\psi)$. We point out that since $U$ and $N(U)$ hold the same linear constraint as in [9], the linear semigroup $e^{-tL}U_0$ and the Duammel term $\int_0^t e^{-(t-s)L}N(U(s)) \mathrm{d}s$ do not include terms which are time-independent or unbounded in $L^1$. Consequently, the $L^1$ estimate $(1.22)$ is obtained from the following integral equation of $U$ 
\[
U(t)=e^{-tL}U(0)+\int_0^t e^{-(t-s)L}N(U(s)) \mathrm{d}s.
\]
The lower $L^1$ estimate $(1.23)$ is obtained by the lower $L^2$ estimate $(1.15)$, the estimate $(1.19)$ with $p=\infty$, and the interpolation inequality.

\vspace{1ex}
We mention that this proof does not need the non-local operator in this reformation. Therefore we remove the difficulty in the case $p=1$, and simplify the analysis of the solutions around the motionless state.  

\vspace{1ex}
This paper is organized as follows. In Section 2 we introduce some notations and function spaces. In Section 3 we state the main result of this paper on the $L^1$  estimates of the solution. In Section 4 we show the $L^1$ estimates.  In Section 5, we derive the $L^1$ estimate of the high frequency part of the linear semigroup. 

\section{Notation}
In this section, we prepare notations and function spaces which will be used throughout the paper. 
$L^p~(1\leq p\leq\infty)$ denotes the usual Lebesgue space on $\mathbb{R}^3$,~and its norm is denoted by $\|\cdot\|_{L^p}$. 
Similarly $W^{m,p} (1\leq p\leq \infty, m\in\{0\}\cup\mathbb{N})$ denotes the $m$-th order $L^p$ Sobolev space on $\mathbb{R}^3$,~and its norm is denoted by $\|\cdot\|_{W^{m,p}}$. We define $H^m=W^{m,2}$ for an integer $m\geq0$. For simplicity, we denote $L^p=L^p\times (L^p)^3\times  (L^p)^9$ (resp. $H^m=H^m\times (H^m)^3\times  (H^m)^9$).

The inner product of $L^2$ is denoted by
$$(f,g):=\int_{\mathbb{R}^3} f(x)\overline{g(x)}dx,~f,g\in L^2.$$
Here the symbol $\overline{\cdot}$ stands for its complex conjugate.
Partial derivatives of a function $u$ in $x_j~(j=1,2,3)$ and $t$ are denoted by $\partial_{x_j}u$ and $\partial_tu$, respectively. $\Delta$ denotes the usual Laplacian with respect to $x$. For a multiindex $\alpha=(\alpha_1,\alpha_2,\alpha_3)\in(\{0\}\cup\mathbb{N})^3$ and $\xi={}^\top(\xi_1,\xi_2,\xi_3)\in\mathbb{R}^3$, we define $\partial_x^\alpha$ and $\xi^\alpha$ as $\partial_x^\alpha =\partial_{x_1}^{\alpha_1} \partial_{x_2}^{\alpha_2} \partial_{x_3}^{\alpha_3}$ and $\xi^\alpha=\xi_1^{\alpha_1}\xi_2^{\alpha_2}\xi_3^{\alpha_3}$, respectively. For a function $u$ and a nonnegative integer $k$, $\nabla^k u$ stands for $\nabla^k u=\{\partial_x^\alpha u|~|\alpha|=k\}$.

For a scalar valued function $\rho=\rho(x)$, we denote by $\nabla\rho$ its gradient with respect to $x$. For a vector valued function $w=w(x)={}^\top(w^1(x),w^2(x),w^3(x))$, we denote by $\mathrm{div}w$ and $(\nabla w)^{jk}=(\partial_{x_k}w^j)$ its divergence and Jacobian matrix with respect to $x$, respectively. For a $3\times3$-matrix valued function $F=F(x)=(F^{jk}(x))$, we define its divergence $\mathrm{div}F$, trace $\mathrm{tr}F$ and determinant $\mathrm{det}F$ by $(\mathrm{div}F)^j=\sum_{k=1}^3 \partial_{x_k}F^{jk}$, $\mathrm{tr}F=\sum_{k=1}^3 F^{kk}$ and $\mathrm{det}F=\sum_{\sigma\in S_3}\mathrm{sgn}(\sigma)F^{1\sigma(1)}F^{2\sigma(2)}F^{3\sigma(3)}$, respectively. Here $S_3$ denotes a third-order symmetric group; for a permutation $\sigma\in S_3$, we denote by $\mathrm{sgn}(\sigma)$ its signature. 

For functions $f=f(x)$ and $g=g(x)$, we denote by $f* g$ the convolution of $f$ and $g$ : 
\[
(f* g)(x)=\int_{\mathbb{R}^3} f(x-y)g(y)\mathrm{d}y.
\]

We denote by $\hat{f}$ or $\mathcal{F}f$ the Fourier transform of a function $f=f(x)$:
\[
\hat{f}(\xi)=(\mathcal{F}f)(\xi)=\frac{1}{(2\pi)^{\frac{3}{2}}}\int_{\mathbb{R}^3} f(x)e^{-i\xi\cdot x} \mathrm{d}x~(\xi\in\mathbb{R}^3).
\] 
The Fourier inverse transform is denoted by $\mathcal{F}^{-1}$:
\[
(\mathcal{F}^{-1}f)(x)=\frac{1}{(2\pi)^{\frac{3}{2}}}\int_{\mathbb{R}^3} f(\xi)e^{i\xi\cdot x} \mathrm{d}\xi~(x\in\mathbb{R}^3).
\] 

We recall the Sobolev inequalities.
\begin{lem}\label{Sobolevineq}
The following inequalities hold:
\begin{align*}
\mathrm{(i)}~\|u\|_{L^p}&\leq C\|u\|_{H^1}~\mathrm{for}~2\leq p\leq6,~u\in H^1.\\[1ex]
\mathrm{(ii)}~\|u\|_{L^p}&\leq C\|u\|_{H^2}~\mathrm{for}~2\leq p \leq\infty,~u\in H^2.
\end{align*}
\end{lem}

We next introduce the following elementary inequality to control the Duammel term. 
\begin{lem}\label{betafunc1}
The following estimates hold:
\[
\mathrm{(i)}~\int_{0}^{t} (1+t-s)^{\frac{1}{2}}(1+s)^{-2}\mathrm{d}s\leq C(1+t)^{\frac{1}{2}},~t\geq0,
\]
Here $C$ is a positive constant independent of $t$.
\end{lem}
\section{Main Result}
In this section, we state the main result of this paper.

We set
$u(t)=(\phi(t),w(t),G(t))=(\rho(t)-1,v(t),F(t)-I).$
Then $u(t)$ satisfies the following initial value problem
\begin{equation}
\left\{\label{system2}
\begin{array}{l}
\partial_t\phi+\mathrm{div} w=g_1(u),\\[1ex]
\partial_t  w-\nu{\Delta} w-\tilde \nu\nabla \mathrm{div}  w+\gamma^2\nabla\phi-\beta^2\mathrm{div} G = g_2(u), \\[1ex]
\partial_t G -\nabla w= g_3(u),\\[1ex]
\nabla \phi +\mathrm{div}{}^\top\! G = g_4(u),\\[1ex]
u|_{t=0}=u_0=(\phi_0,w_0,G_0).
\end{array}
\right.
\end{equation}
Here $g_j(u), j=1, 2, 3, 4,$ denote the nonlinear terms;
\begin{gather*}
\begin{array}{l}
g_1(u)=-\mathrm{div}(\phi w), \\[1ex]
\displaystyle
g_2(u)=-w\cdot\nabla w+\frac{\phi}{1+\phi}(-\nu\Delta w-\tilde{\nu}\nabla\mathrm{div}w+\gamma^2\nabla \phi)-\frac{1}{1+\phi}\nabla R(\phi)\\[2ex]
\displaystyle
\quad\quad-\frac{\beta^2\phi}{1+\phi}\mathrm{div} G +\frac{\beta^2}{1+\phi}\mathrm{div}(\phi G+G{}^\top\! G+\phi G{}^\top\! G), \\[2ex]
g_3(u)=-w\cdot\nabla G+\nabla w G, \\[1ex]
g_4(u)=g_4(\phi,G)=-\mathrm{div}(\phi{}^\top\! G),\\[1ex] 
\end{array}   
\end{gather*}
where
\[
R(\phi)=\phi^2\int_0^1 {P}^{\prime\prime}(1+s\phi)\mathrm{d}s,~\nabla R=O(\phi)\nabla \phi
\]
for $|\phi|\ll1$. 

The following proposition ensures the global in time existence and the $L^2$ decay estimates of solutions.

\begin{prop}\label{GE}{\sc $([6,12])$}
Let $u_0\in H^N,~N\geq3$. There is a positive number $\epsilon_0$ such that if $u_0$ satisfies $\|u_0\|_{L^1}+\|u_0\|_{H^3}\leq\epsilon_0$, then there exists a unique solution $u(t)\in C([0,\infty);H^N)$ of the problem $(3.1)$, and $u(t)=(\phi(t),w(t),G(t))$ satisfies 
\begin{gather*}
\| u(t)\|_{H^N}^2+\int_{0}^{t}(\|\nabla\phi(s)\|_{H^{N-1}}^2+\|\nabla w(s)\|_{H^N}^2+\|\nabla G(s)\|_{H^{N-1}}^2)\mathrm{d}s\leq C_N\|u_0\|_{H^N}^2,
\\
\|\nabla^k u(t)\|_{L^2}\leq C(1+t)^{-\frac{3}{4}-\frac{k}{2}}(\|u_0\|_{L^1}+\|u_0\|_{H^N}) 
\end{gather*}
for $k=0,1,2,\ldots,N-1$ and $t\geq 0$.

In addition, if there exists a positive number $r>0$ such that the following condition satisfies
\begin{align}
|\hat{\phi}_0(\xi)|>c_0, |\hat{m}_0(\xi)|+|\hat{\mathcal{G}}_0(\xi)-{}^\top\hat{\mathcal{G}}_0(\xi)|\leq  c_1|\xi|^{\eta_0} \label{lowfreqGE}
\end{align}
for $0\leq |\xi|\leq r$, where $(m_0,\mathcal{G}_0):=(\rho_0 v_0, \rho_0 F_0-I)$; $c_0$, $c_1$ and $\eta_0$ are positive numbers independent of $t$, the following lower $L^2$ estimate holds
\begin{align}
\|u(t)\|_{L^2}\geq c(1+t)^{-\frac{3}{4}} \label{L2sharpestimate}
\end{align}
uniformly for $t\geq R_1$. Here $R_1$ is a large positive number, and $c$ is a positive number independent of $t$.
\end{prop}

We next state the main result of this paper which gives the upper estimate of the $L^1$ norm of solutions.

\begin{thm}\label{mainthm}

Assume that $\phi_0$ and $G_0$ satisfy $\nabla \phi_0 +\mathrm{div}{}^\top\! G_0 =g_4(\phi_0,G_0)$ and $(I+G_0)^{-1}=\nabla X_0$ for some vector field $X_0$. There is a positive number $\epsilon$ such that if $u_0=(\phi_0,w_0,G_0)$ satisfies $\|u_0\|_{H^4}+\|u_0\|_{L^1}\leq\epsilon$ and $u_0\in W^{2,1}$, then there exists a unique solution $u(t)\in C([0,\infty);H^4)$ of the problem $(3.1)$ satisfying
\begin{align*}
\|u(t)\|_{L^1}\leq C(1+t)^{\frac{1}{2}}(\|u_0\|_{W^{2,1}}+\|u_0\|_{H^4}) 
\end{align*}
uniformly for $t>0$.
Here $C$ is a positive constant.
\end{thm}

We also prove the lower $L^1$ estimate to show that the obtained rate $(1+t)^{\frac{1}{2}}$ is sharp.

\begin{thm}\label{sharpest}

Under the assumptions in Proposition 3.1 and Theorem 3.2, there exists a positive large number $R_1$ such that the lower $L^1$ estimate of the solution $u$ of the problem $(3.1)$ holds
\begin{align*}
\|u(t)\|_{L^1}\geq c(1+t)^{\frac{1}{2}} 
\end{align*}
uniformly for $t\geq R_1$.
Here $c$ is a positive constant independent of time $t$.
\end{thm}

{\bf Proof of Theorem 3.3.} 
Since we assume the condition $(3.2)$, the lower $L^2$ estimate $(3.3)$ holds uniformly for $t\geq R_1$. Here $R_1$ is taken in Proposition 3.1. We next introduce the following $L^\infty$ decay estimate which can be proved in a similar manner to [9, Theorem 3.2 (i)] and Theorem 3.2  
\begin{align}
\|u(t)\|_{L^\infty}\leq C(1+t)^{-2}(\|u_0\|_{L^1}+\|u_0\|_{H^3}),~t\geq0. \label{Linftyestimate}
\end{align}
By using $(3.3)$,  $(3.4)$ and the interpolation inequality $\|u(t)\|_{L^2}\leq  \|u(t)\|_{L^1}^{\frac{1}{2}}\|u(t)\|_{L^\infty}^{\frac{1}{2}}$, we have 
\begin{align*}
c(1+t)^{-\frac{3}{4}}&\leq \|u(t)\|_{L^2} \\
&\leq  \|u(t)\|_{L^\infty}^{\frac{1}{2}}\|u(t)\|_{L^1}^{\frac{1}{2}} \\
&\leq  C(1+t)^{-1} \|u(t)\|_{L^1}^{\frac{1}{2}}. 
\end{align*}
This yields
\[
\|u(t)\|_{L^1}\geq c(1+t)^{\frac{1}{2}}. 
\]
This complete the proof. $\blacksquare$

\section{Proof of Theorem 3.2}

This section is devoted to the proof of Theorem 3.2.

The global existence of solutions is proved by Proposition 3.1. Hence we focus on the derivation of the $L^1$ estimate. 

As we mentioned above, it is not valid to apply the semigroup theory to the problem $(3.1)$ directly. To overcome this difficulty, we reformulate the problem $(3.1)$.

Let $x=x(X,t)$ be a material coordinate which solves the following flow map
\[
\left\{
\begin{array}{l}
\displaystyle
\frac{dx}{dt}(X,t)=v(x(X,t),t),\\[1ex]
x(X,0)=X.
\end{array}
\right.
\]
We next define $X=X(x,t)$ as the inverse of $x=x(X,t)$ and $\psi=x-X$.
According to [1,21], $F$ is defined as $F=\frac{\partial x}{\partial X}$. It is shown in [19] that its inverse $F^{-1}$ is written as $F^{-1}(x,t)=\nabla X(x,t)$ if $F_0^{-1}$ has the form $F^{-1}_0=\nabla X_0$ with some vector field $X_0$.  
Then $\psi$ solves
\begin{align}
\partial_t\psi-v=-v\cdot\nabla\psi, \label{eqpsi}
\end{align}
and satisfies
\begin{equation}
G=\nabla\psi+h(\nabla\psi),\label{3.1}
\end{equation}
where $h(\nabla\psi)=(I-\nabla\psi)^{-1}-I-\nabla\psi$. 

We note that $(4.2)$ is equivalent to
\begin{equation}
\nabla\psi=I-(I+G)^{-1}.\label{3.2}
\end{equation}

We next read the constraint $(1.8)$ as 
\[
1+\phi=\rho=\mathrm{det}F^{-1}=\mathrm{det}(I+G)^{-1}.
\]
By using $(4.3)$ and the following expansion for a $3\times3$ matrix $A$
\[
\mathrm{det}(I+A)=1+\mathrm{tr}A+\frac{1}{2}(\mathrm{tr}A^2-(\mathrm{tr}A)^2)+\mathrm{det}A,
\]
we have
\begin{align}
\phi=-\mathrm{div}\psi+\frac{1}{2}(\mathrm{tr}(\nabla\psi)^2-(\mathrm{tr}(\nabla\psi))^2)-\mathrm{det}(\nabla\psi). \label{3.3}
\end{align}

We set $\psi_0=\psi|_{t=0}$. The following estimates hold for $\phi$, $G$ and $\nabla\psi$.
\begin{lem}\label{equivforGpsi}
Assume that $G$ and $\psi$ satisfy $(4.2)$. There is a positive number $\epsilon_1<1$ such that if $\|G\|_{C([0,\infty);H^3)}\leq\epsilon_1$, the following inequalities hold uniformly for $t\geq 0$:
\begin{align}
&C^{-1}\|\nabla\psi(t)\|_{L^p}\leq\|G(t)\|_{L^p} \leq C\|\nabla\psi(t)\|_{L^p},~p=1,2,  \label{equiv2} \\
&\| \nabla^2\psi(t)\|_{L^2} \leq C\|\nabla G(t)\|_{L^2}, \label{equiv1.1}  \\
&\| \nabla^3\psi(t)\|_{L^2} \leq C(\|\nabla G(t)\|_{H^1}^2+\|\nabla^2 G(t)\|_{L^2}), \label{equiv1.2}  \\
&\| \nabla^4\psi(t)\|_{L^2} \leq C(\|\nabla G(t)\|_{H^1}\|\nabla^2 G(t)\|_{H^1}+\|\nabla^3 G(t)\|_{L^2}). \label{equiv1.3}  \\
&\| \nabla^2\psi_0\|_{L^1} \leq C\|\nabla G_0\|_{L^1}, \label{equiv1.1.1}  \\
&\| \nabla^3\psi_0\|_{L^1} \leq C\|\nabla^2 G_0\|_{L^1}, \label{equiv1.1.2}  \\
&\| \phi(t)\|_{L^1}\leq \|\mathrm{div}\psi(t)\|_{L^1}+C(1+\|\nabla\psi(t)\|_{L^\infty})\|\nabla\psi(t)\|_{L^2}^2. \label{equiv1.4}
\end{align}
\end{lem}
{\bf Proof.} The inequalities $(4.5)$--$(4.8)$ are shown in [9, Lemma 4.1]. The inequalities $(4.9)$ and $(4.10)$ are established in a similar argument as in the proof of \eqref{equiv1.1}. We only derive $(4.11)$ here.

In order to obtain $(4.11)$, we use $(4.4)$. Since
\[
\|(\mathrm{tr}(\nabla\psi)^2-(\mathrm{tr}(\nabla\psi))^2)\|_{L^1}\leq C\|\nabla\psi\|_{L^2}^2,
\]
\[
\|\mathrm{det}(\nabla\psi)\|_{L^1}\leq C\|\nabla\psi\|_{L^\infty}\|\nabla\psi\|_{L^2}^2,
\]
we have
\[
\| \phi\|_{L^1}\leq \|\mathrm{div}\psi\|_{L^1}+C(1+\|\nabla\psi\|_{L^\infty})\|\nabla\psi\|_{L^2}^2.
\]
This completes the proof of Lemma 4.1. $\blacksquare$
\vspace{1ex}

Let $U$ and $U_0$ be $U=(\tilde{\phi},w,\tilde{G})=(-\mathrm{div}\psi,w,\nabla\psi)$ and $U_0=(\tilde{\phi}_0,w_0,\tilde{G}_0)=(-\mathrm{div}\psi_0,w_0,\nabla\psi_0)$, respectively. By taking $\epsilon_0$ in Proposition 3.1 such that $C_3\epsilon_0\leq \epsilon_1$, thanks to Lemma 4.1, the behavior of $u$ is identified from $U$: $\|u(t)\|_{L^1}\leq C\|U(t)\|_{L^1}$. By coupling $(3.1)$ and $(4.1)$, we arrive at the following problem for $U$
\begin{equation}
\left\{\label{system3}
\begin{array}{l}
\partial_t\tilde{\phi}+\mathrm{div} w=N_1(U),\\
\partial_t  w-\nu{\Delta} w-\tilde \nu\nabla \mathrm{div}  w+\gamma^2\nabla\tilde{\phi}-\beta^2\mathrm{div}\tilde{G} = N_2(U), \\
\partial_t \tilde{G} - \nabla w=  N_3(U),\\
\tilde{\phi} +\mathrm{tr}\tilde{G} =0,~\tilde{G}=\nabla\psi,\\
U|_{t=0}=U_0.
\end{array}
\right.
\end{equation}
Here $N_j(U), j=1, 2, 3,$ denote nonlinear terms;
\begin{gather*}
\begin{array}{l}
N_1(U)=\mathrm{div}(w\cdot\nabla\psi), \\[1ex]
\displaystyle
N_2(U)=g_2(u)-\frac{\gamma^2}{2}\nabla(\mathrm{tr}(\nabla\psi)^2-(\mathrm{tr}(\nabla\psi))^2)-\gamma^2\nabla\mathrm{det}(\nabla\psi)+\beta^2\mathrm{div}h(\nabla\psi),\\[2ex]
N_3(U)=-\nabla(w\cdot\nabla\psi).
\end{array}   
\end{gather*}
We note that $N_1(U)$ and $N_3(U)$ also satisfy the same linear constraint as $\tilde{\phi}$ and $\tilde{G}$
\begin{align}
N_1(U) +\mathrm{tr}N_3(U) =0.  
\end{align}

In what follows, we omit tildes $\tilde{\cdot}$ of $\tilde{\phi}$ and $\tilde{G}$ for simplicity. The problem $(4.12)$ is reduced to 
\begin{equation}
\left\{\label{problem5}
\begin{array}{l}
\partial_t U+LU=N(U), \\
\phi +\mathrm{tr}G = 0,~G=\nabla\psi,\\
U|_{t=0}=U_0,
\end{array}
\right.
\end{equation}
where
$$
L=
\left(
\begin{array}{@{\ }cc@{\ }cc@{\ } }
0 & \mathrm{div}& 0\\
\gamma^2\nabla & -\nu\Delta-\tilde{\nu}\nabla\mathrm{div} & -\beta^2\mathrm{div} \\
0 & -\nabla & 0 
\end{array}
\right),~
N(U)=\left(
\begin{array}{l}
N_1(U) \\
N_2(U) \\
N_3(U)
\end{array}
\right).
$$
By using the Duammel principle, $U$ satisfies the following integral equations
\begin{align}
U(t)=e^{-tL}U_0+\int_0^t e^{-(t-s)L}N(U(s)) \mathrm{d}s. \label{IECh4}
\end{align}

To analyize the linearized semigroup $e^{-tL}$, we introduce several notations. We set
\begin{align*}
Q&=I-P=\mathcal{F}^{-1}\frac{\xi{}^\top\xi}{|\xi|^2}\mathcal{F}, \\
\mathcal{K}_t^\lambda(x)&=\mathcal{F}^{-1}\left[\frac{e^{\lambda_+(\xi)t}-e^{\lambda_-(\xi)t}}{\lambda_+(\xi)-\lambda_-(\xi)}\right](x),\\
\mathcal{K}_t^\mu(x)&=\mathcal{F}^{-1}\left[\frac{e^{\mu_+(\xi)t}-e^{\mu_-(\xi)t}}{\mu_+(\xi)-\mu_-(\xi)}\right](x).
\end{align*}
Here $\lambda_{\pm}(\xi)$ and $\mu_{\pm}(\xi)$ are given by
\begin{align*}
\lambda_{\pm}(\xi)&=\frac{-\nu|\xi|^2\pm\sqrt{\nu^2|\xi|^4-4\beta^2|\xi|^2} }{2},\\
\mu_{\pm}(\xi)&=\frac{-(\nu+\tilde{\nu})|\xi|^2\pm\sqrt{(\nu+\tilde{\nu})^2|\xi|^4-4(\beta^2+\gamma^2)|\xi|^2} }{2}.
\end{align*} 

We see that the following properties of $\lambda_{\pm}(\xi)$ and $\mu_{\pm}(\xi)$ hold:
\begin{align*}
&\lambda_{\pm}(\xi)\sim -\frac{\nu}{2}|\xi|^2\pm i\beta|\xi|,~\mathrm{for}~|\xi|\ll1, \\
&\lambda_+(\xi)\sim-\frac{\beta^2}{\nu},~\lambda_-(\xi)\sim-\nu|\xi|^2,~\mathrm{for}~|\xi|\gg1,\\
&\mu_\pm(\xi)\sim-\frac{\nu+\tilde{\nu}}{2}|\xi|^2\pm i\sqrt{\beta^2+\gamma^2}|\xi|,~\mathrm{for}~|\xi|\ll1, \\
&\mu_+(\xi)\sim-\frac{\beta^2+\gamma^2}{\nu+\tilde{\nu}},~\mu_-(\xi)\sim-(\nu+\tilde{\nu})|\xi|^2,~\mathrm{for}~|\xi|\gg1. 
\end{align*}

According to [9], the expression of the semigroup $U(t)={}^\top(\phi(t),w(t),G(t))=e^{-tL}U_0$ is given by
\begin{align}
\left(
\begin{array}{l}
\phi(x,t) \\
w(x,t) \\
G(x,t)
\end{array}
\right)=
\left(
\begin{array}{@{\ }cc@{\ }cc@{\ } }
\mathcal{K}^{11}(t) & \mathcal{K}^{12}(t)& \mathcal{K}^{13}(t)\\
\mathcal{K}^{21}(t) & \mathcal{K}^{22}(t)& \mathcal{K}^{23}(t)\\
\mathcal{K}^{31}(t) & \mathcal{K}^{32}(t)& \mathcal{K}^{33}(t)
\end{array}
\right)
\left(
\begin{array}{l}
\phi_0(x) \\
w_0(x) \\
G_0(x)
\end{array}
\right),
\label{solformula}
\end{align}
provided that $\phi_0+\mathrm{tr}G_0=0$ and $G_0=\nabla\psi_0$.

Here $\mathcal{K}^{jk}(t), ~j,k=1,2,3,$ are the linear operators defined as
\begin{align*}
\mathcal{K}^{11}(t)\phi_0(x)&=(\partial_t -(\nu+\tilde\nu)\Delta)( \mathcal{K}_t^\mu*\phi_0)(x)\\[0.5ex]
&=\mathcal{F}^{-1}\left[\frac{\mu_+(\xi)e^{\mu_-(\xi)t}-\mu_-(\xi)e^{\mu_+(\xi)t}}{\mu_+(\xi)-\mu_-(\xi)}\hat{\phi}_0(\xi)\right](x),\\[0.5ex]
\mathcal{K}^{12}(t)w_0(x)&=-\mathrm{div}(\mathcal{K}_t^\mu*w_0)(x) \\[0.5ex]
&=i\mathcal{F}^{-1}\left[\frac{e^{\mu_+(\xi)t}-e^{\mu_-(\xi)t}}{\mu_+(\xi)-\mu_-(\xi)}\xi\cdot\hat{w}_0(\xi)\right](x), \\[0.5ex] 
\mathcal{K}^{13}(t)G_0(x)&=0, \\[0.5ex]
\mathcal{K}^{21}(t)\phi_0(x)&=-\gamma^2\nabla(\mathcal{K}_t^\mu*\phi_0)(x), \\[0.5ex]
&=i\gamma^2\mathcal{F}^{-1}\left[\frac{e^{\mu_+(\xi)t}-e^{\mu_-(\xi)t}}{\mu_+(\xi)-\mu_-(\xi)}\hat{\phi}_0(\xi)\xi\right](x), \\[0.5ex] 
\mathcal{K}^{22}(t)w_0(x)&=\partial_t (\mathcal{K}_t^\lambda*Pw_0)(x)+\partial_t (\mathcal{K}_t^\mu*Qw_0)(x)\\[0.5ex]
&=\mathcal{F}^{-1}\left[\frac{\lambda_+(\xi)e^{\lambda_+(\xi)t}-\lambda_-(\xi)e^{\lambda_-(\xi)t}}{\lambda_+(\xi)-\lambda_-(\xi)}\left(I-\frac{\xi^{\top}\xi}{|\xi|^2}\right)\hat{w}_0(\xi)\right](x)  \\[0.5ex]
&\quad+\mathcal{F}^{-1}\left[\frac{\mu_+(\xi)e^{\mu_+(\xi)t}-\mu_-(\xi)e^{\mu_-(\xi)t}}{\mu_+(\xi)-\mu_-(\xi)}\frac{\xi^{\top}\xi}{|\xi|^2}\hat{w}_0(\xi)\right](x), \\[0.5ex]
\mathcal{K}^{23}(t)G_0(x)&=\beta^2(PG_0*\nabla \mathcal{K}_t^\lambda)(x)+\beta^2(QG_0*\nabla \mathcal{K}_t^\mu)(x) \\[0.5ex]
&=-i\beta^2\mathcal{F}^{-1}\left[\frac{e^{\lambda_+(\xi)t}-e^{\lambda_-(\xi)t}}{\lambda_+(\xi)-\lambda_-(\xi)}\left(I-\frac{\xi^{\top}\xi}{|\xi|^2}\right)\hat{G}_0(\xi)\xi\right](x)  \\[0.5ex]
&\quad\quad-i\beta^2\mathcal{F}^{-1}\left[\frac{e^{\mu_+(\xi)t}-e^{\mu_-(\xi)t}}{\mu_+(\xi)-\mu_-(\xi)}\frac{\xi^{\top}\xi}{|\xi|^2}\hat{G}_0(\xi)\xi\right](x),\\[0.5ex]
\mathcal{K}^{31}(t)\phi_0(x)&=0,\\[0.5ex]
\mathcal{K}^{32}(t)w_0(x)&=(Pw_0*\nabla \mathcal{K}_t^\lambda)(x)+(Qw_0*\nabla \mathcal{K}_t^\mu)(x) \\[0.5ex]
&=-i\mathcal{F}^{-1}\left[\frac{e^{\lambda_+(\xi)t}-e^{\lambda_-(\xi)t}}{\lambda_+(\xi)-\lambda_-(\xi)}\left(I-\frac{\xi^{\top}\xi}{|\xi|^2}\right)\hat{w}_0(\xi){}^\top\xi\right](x)  \\[0.5ex]
&\quad\quad-i\mathcal{F}^{-1}\left[\frac{e^{\mu_+(\xi)t}-e^{\mu_-(\xi)t}}{\mu_+(\xi)-\mu_-(\xi)}\frac{\xi^{\top}\xi}{|\xi|^2}\hat{w}_0(\xi){}^\top\xi\right](x),\\[2ex]
\mathcal{K}^{33}(t)G_0(x)&=(\partial_t -\nu\Delta )(\mathcal{K}_t^\lambda*PG_0)(x)\\[0.5ex]
&\quad +(\partial_t -(\nu+\tilde\nu)\Delta)( \mathcal{K}_t^\mu*QG_0)(x)\\[0.5ex]
&=\mathcal{F}^{-1}\left[\frac{\lambda_+(\xi)e^{\lambda_-(\xi)t}-\lambda_-(\xi)e^{\lambda_+(\xi)t}}{\lambda_+(\xi)-\lambda_-(\xi)}\left(I-\frac{\xi^{\top}\xi}{|\xi|^2}\right)\hat{G}_0(\xi)\right](x)  \\[0.5ex]
&\quad\quad+\mathcal{F}^{-1}\left[\frac{\mu_+(\xi)e^{\mu_-(\xi)t}-\mu_-(\xi)e^{\mu_+(\xi)t}}{\mu_+(\xi)-\mu_-(\xi)}\frac{\xi^{\top}\xi}{|\xi|^2}\hat{G}_0(\xi)\right](x).
\end{align*}

In view of the asymptotic profiles of $\lambda_{\pm}(\xi)$ and $\mu_{\pm}(\xi)$ in the Fourier space, we decompose the solution $U(t)$ of the problem $(4.14)$ into its low and high frequency parts. Let $\hat\varphi_1,~\hat\varphi_M,~ \hat\varphi_\infty\in C^\infty(\mathbb{R}^3;[0,1])$ be cut-off functions such that
\begin{align*}
\hat\varphi_1(\xi)&=
\begin{cases}
1 & |\xi|\leq \frac{M_1}{2}, \\
0 & |\xi|\geq \frac{M_1}{\sqrt{2}},
\end{cases}
~\hat\varphi_1(-\xi)=\hat\varphi_1(\xi),\\
\hat\varphi_\infty(\xi)&=
\begin{cases}
1 & |\xi|\geq 2M_2 \\
0 & |\xi|\leq \sqrt{2}M_2,
\end{cases}
~\hat\varphi_\infty(-\xi)=\hat\varphi_\infty(\xi), \\
\hat\varphi_M(\xi)&=1-\hat\varphi_1(\xi)-\hat\varphi_\infty(\xi),
\end{align*}
where
\[
M_1=\min\left\{\frac{\beta}{\nu},\frac{\sqrt{\beta^2+\gamma^2}}{\nu+\tilde\nu}\right\},~
M_2=\max\left\{\frac{\beta}{\nu},\frac{\sqrt{\beta^2+\gamma^2}}{\nu+\tilde\nu}\right\}.
\]
We define operators $P_j,~j=1,\infty,$ on $L^2$ as
\[
P_1u=\mathcal{F}^{-1}(\hat\varphi_1\hat{u}),~P_\infty u=\mathcal{F}^{-1}((\hat\varphi_M+\hat\varphi_\infty)\hat{u})~\mathrm{for}~u\in L^2.
\]

\begin{lem}\label{lemP1Pinfty}
$P_j~(j=1,\infty)$ have the following properties.

$(i)$ $P_1+P_\infty=I$.

$(ii)$ $\partial_x^\alpha P_1=P_1\partial_x^\alpha$, $\|\partial_x^\alpha P_1 f\|_{L^2}\leq C_\alpha \|f\|_{L^2}$ for $\alpha\in(\{0\}\cup\mathbb{N})^3$ and $f\in L^2$.

$(iii)$ $\partial_x^\alpha P_\infty=P_\infty\partial_x^\alpha$, $\|\partial_x^\alpha P_\infty f\|_{L^2}\leq C \|\nabla\partial_x^\alpha P_\infty f\|_{L^2}$ for $\alpha\in(\{0\}\cup\mathbb{N})^3$ with $|\alpha|=k\geq0$ and $f\in H^{k+1}$.
\end{lem}
Lemma 4.2 immediately follows from the definitions of $P_j,~j=1,\infty,$ and the Plancherel theorem. We omit the proof. 
 
The solution $U(t)$ of $(4.14)$ is decomposed as
\[
U(t)=U_1(t)+U_\infty(t),~U_1(t)=P_1U(t),~U_\infty(t)=P_\infty U(t).
\]
By applying $P_j$ to $(4.15)$, $U_j(t)=(\phi_j(t),w_j(t),G_j(t)),~j=1,\infty,$ satisfy
\begin{gather}
\left\{
\begin{array}{l}
\displaystyle
U_j(t)=e^{-tL}U_j(0)+\int_{0}^{t} e^{-(t-s)L}P_jN(U(s))\mathrm{d}s, \\[1.5ex]
\phi_j+\mathrm{tr}G_j=0,~P_jN_1(U)+\mathrm{tr}P_jN_3=0, \\[1.5ex]
U_j|_{t=0}=P_jU_0.
\end{array}
\right.
\label{jfreq}  
\end{gather}

Concerning to the $L^1$ estimate of $e^{-tL}U_0$, we have the following proposition.

\begin{prop}\label{L1estetL}
Let $\phi_0+\mathrm{tr}G_0=0$ and $G_0=\nabla\psi_0$. Then the following estimates hold for $t>0$:
\begin{align*}
(i) \hspace{3ex} \|e^{-tL}P_1U_0\|_{L^1}&\leq C(1+t)^{\frac{1}{2}}\|U_0\|_{L^1}, \\
(ii)~\|e^{-tL}P_\infty U_0\|_{L^1}&\leq Ce^{-ct}\|U_0\|_{W^{2,1}}. 
\end{align*}
\end{prop}

The estimate $(i)$ is done in [9, Lemma 6.4]. For the estimate $(ii)$, we will discuss in Section 5.

\vspace{2ex}
In order to estimate the Duammel terms $\int_{0}^{t} e^{-(t-s)L}P_jN(U(s))\mathrm{d}s,~j=1,\infty$, we give the following $L^2$ decay estimates for $\nabla^k U(t)$.
\begin{prop}\label{H2estimate}
There exists a positive number $\epsilon_1$ such that if $\|u_0\|_{L^1}+\|u_0\|_{H^4}\leq\epsilon_1$, then the following inequality holds for $k=0,1,2,3$ and $t\geq 0$:
$$
\|\nabla^k U(t)\|_{L^2}\leq C(1+t)^{-\frac{3}{4}-\frac{k}{2}}(\|u_0\|_{L^1}+\|u_0\|_{H^4}),
$$
\end{prop}
Proposition 4.4 follows from Proposition 3.1 and Lemma 4.1.

The estimates of $\int_{0}^{t} \|e^{-(t-s)L}P_jN(U(s))\|_{L^1}\mathrm{d}s, j=1,\infty,$ are given as follows.

\begin{lem}\label{estetLPjN}
There exists a positive number $\epsilon_1$ such that if $u_0\in L^1\cap H^4$ and $\|u_0\|_{H^4}\leq\epsilon_1$, then the following estimates hold:
\[
\int_{0}^{t} \|e^{-(t-s)L}P_jN(U(s))\|_{L^1}\mathrm{d}s\leq C(1+t)^{\frac{1}{2}}(\|u_0\|_{L^1}+\|u_0\|_{H^4}),~j=1,\infty,~t\geq0.
\]
\end{lem}

{\bf Proof.} We first show the case $j=1$. We obtain the following estimate in a similar argument as in the proof of Lemma 4.3: 
\begin{align}
\begin{array}{l}
\displaystyle
\left\|e^{-(t-s)L}P_1N(U(s))\right\|_{L^1}
\leq C(1+t-s)^{\frac{1}{2}}\|N(U(s))\|_{L^1} .
\end{array}
\label{estN111L} 
\end{align}
In view of Lemma 2.1 and Proposition 4.4, we have
\begin{align}
\begin{array}{l}
\|N(U(s))\|_{L^1}\leq C\|U(s)\|_{H^2}\|\nabla U(s)\|_{H^1}
\leq C(1+s)^{-2}(\|u_0\|_{L^1}+\|u_0\|_{H^4}).
\end{array}
\label{ineqLEM6}
\end{align}
We see from Lemma 2.2, $(4.18)$ and $(4.19)$ that
\begin{align}
\begin{array}{l}
\displaystyle
\int_0^t \|e^{-(t-s)L}P_1N(U(s))\|_{L^1}\mathrm{d}s \\[2ex]
\quad
\displaystyle
\leq
\displaystyle
C\int_0^t (1+t-s)^{\frac{1}{2}}\|N(U(s))\|_{L^1}\mathrm{d}s \\[2ex]
\quad
\displaystyle
\leq
C\int_0^t (1+t-s)^{\frac{1}{2}}(1+s)^{-2}\mathrm{d}s(\|u_0\|_{L^1}+\|u_0\|_{H^4}) \\[2ex]
\quad
\displaystyle
\leq C(1+t)^{\frac{1}{2}}(\|u_0\|_{L^1}+\|u_0\|_{H^4}). 
\end{array}
\label{ineqLEM671}
\end{align}

We next consider the case $j=\infty$.

We obtain the following estimate in a similar argument as in the proof of Lemma 4.3: 
\begin{align}
\|e^{-(t-s)L}P_\infty N(U(s))\|_{L^1}\leq Ce^{-c(t-s)}\|N(U(s))\|_{W^{1,1}}.
\label{ineqLEM6110} 
\end{align}
In view of Lemma 2.1 and Proposition 4.4, we have
\begin{align}
\begin{array}{l}
\|N(U(s))\|_{W^{2,1}}\leq C(\|U(s)\|_{H^3}^2+\|u(s)\|_{L^2}\|\nabla^4 w(s)\|_{L^2})
\leq C\|u_0\|_{H^4}.
\end{array}
\label{ineqLEM6111}
\end{align}
It follows from $(4.21)$ and $(4.22)$ that
\begin{align*}
\begin{array}{l}
\displaystyle
\int_0^t \|e^{-(t-s)L}P_\infty N(U(s))\|_{L^1}\mathrm{d}s \\[2ex]
\quad
\displaystyle
\leq
C\int_0^t e^{-c(t-s)}\|N(U(s))\|_{W^{2,1}}\mathrm{d}s
\\[2ex]
\quad
\displaystyle
\leq
C\|u_0\|_{H^4}
\\[2ex]
\quad
\displaystyle
\leq C(1+t)^{\frac{1}{2}}\|u_0\|_{H^4}.
\end{array}
\end{align*}
This completes the proof. $\blacksquare$

\textbf{Proof of Theorem 3.2.}
By taking $L^1$ norm of the first equation of $(4.17)$, we have
\begin{align}\label{lowfreqineqL1}
\|U_j(t)\|_{L^1}&\leq \|e^{-tL}U_j(0)\|_{L^1}+\int_{0}^{t} \|e^{-(t-s)L}P_j N(s)\|_{L^1}\mathrm{d}s. 
\end{align}
Combining Lemma 4.3, Lemma 4.5 and $(4.23)$, we arrive at
\[
\|U(t)\|_{L^1}\leq C(1+t)^{\frac{1}{2}}(\|u_0\|_{W^{2,1}}+\|u_0\|_{H^4}),~t\geq0.
\]
By using Lemma 4.1, we obtain
\[
\|u(t)\|_{L^1}\leq C\|U(t)\|_{L^1}\leq C(1+t)^{\frac{1}{2}}(\|u_0\|_{W^{2,1}}+\|u_0\|_{H^4}),~t\geq0.
\]
This completes the proof of Theorem 3.2. $\blacksquare$ \\[1ex]

\section{Proof of Proposition 4.3 (ii).} 

In this section, we prove Proposition 4.3 (ii).

For $j=M,\infty$, we set
\begin{align*}
\mathcal{K}_j^{\lambda_{\pm}}(t)f(x)&=\mathcal{F}^{-1}\left[\frac{e^{\lambda_{\pm}(\xi)t}}{\lambda_{+}(\xi)-\lambda_{-}(\xi)}\hat\varphi_j(\xi)\hat{f}(\xi)\right],\\
\mathcal{M}_j^{\lambda_{\pm}}(t)f(x)&=\mathcal{F}^{-1}\left[\frac{\lambda_{\mp}(\xi)e^{\lambda_{\pm}(\xi)t}}{\lambda_{+}(\xi)-\lambda_{-}(\xi)}\hat\varphi_j(\xi)\hat{f}(\xi)\right],\\
\mathcal{K}_j^{\mu_{\pm}}(t)f(x)&=\mathcal{F}^{-1}\left[\frac{e^{\mu_\pm(\xi)t}}{\mu_+(\xi)-\mu_-(\xi)}\hat\varphi_j(\xi)\hat{f}(\xi)\right],\\
\mathcal{M}_j^{\mu_{\pm}}(t)f(x)&=\mathcal{F}^{-1}\left[\frac{\mu_\mp(\xi)e^{\mu_\pm(\xi)t}}{\mu_+(\xi)-\mu_-(\xi)}\hat\varphi_j(\xi)\hat{f}(\xi)\right].
\end{align*}

We first consider the high frequency part.

\begin{lem}\label{highfreqestetL}
The following estimates hold for $\alpha\in(\mathbb{N}\cup\{0\})^3$, $j\geq k\geq0$ and $t>0$:
\begin{align}
&(i)~\left\|\partial_t^j\partial_x^\alpha[\mathcal{K}_{\infty}^{\lambda_+}(t)f]\right\|_{L^1}+\left\|\partial_t^j\partial_x^\alpha[\mathcal{K}_{\infty}^{\mu_+}(t)f]\right\|_{L^1} \leq  Ce^{-ct}\|f\|_{W^{(|\alpha|-1)^+,1}}, \label{Sineq1} 
\\[2ex] 
&(ii)~
\begin{array}{l}
\displaystyle
\left\|\partial_t^j\partial_x^\alpha[\mathcal{K}_{\infty}^{\lambda_-}(t)f]\right\|_{L^1}+\left\|\partial_t^j\partial_x^\alpha[\mathcal{K}_{\infty}^{\mu_-}(t)f]\right\|_{L^1}\\[1ex]
\displaystyle
\quad
\leq Ce^{-ct}t^{-(j-k)}\|f\|_{W^{2k+(|\alpha|-1)^+,1}},
\end{array}
\label{Sineq2} 
\\[2ex] 
&(iii)~\left\|\partial_t^j\partial_x^\alpha[\mathcal{M}_{\infty}^{\lambda_+}(t)f]\right\|_{L^1}+\left\|\partial_t^j\partial_x^\alpha[\mathcal{M}_{\infty}^{\mu_+}(t)f]\right\|_{L^1} \leq Ce^{-ct}\|f\|_{W^{|\alpha|,1}}, \label{Sineq3} 
\\[2ex]  
&(iv)~
\begin{array}{l}
\displaystyle
\left\|\partial_t^j\partial_x^\alpha[\mathcal{M}_{\infty}^{\lambda_-}(t)f]\right\|_{L^1}+\left\|\partial_t^j\partial_x^\alpha[\mathcal{M}_{\infty}^{\mu_-}(t)f]\right\|_{L^1}\\[1ex]
\displaystyle
\quad
\leq Ce^{-ct}t^{-(j-k)}\|f\|_{W^{2k+(|\alpha|-1)^+,1}},
\end{array}  \label{Sineq4} 
\\[2ex]  
&(v)~\left\|\partial_t^j\partial_x^\alpha[\mathcal{K}_{\infty}^{\lambda_+}(t)Qf]\right\|_{L^1}+\left\|\partial_x^\alpha[\mathcal{K}_{\infty}^{\mu_+}(t)Qf]\right\|_{L^1} \leq Ce^{-ct}\|f\|_{W^{|\alpha|,1}},  \label{Sineq5} 
\\[2ex]  
&(vi)~
\begin{array}{l}
\displaystyle
\left\|\partial_t^j\partial_x^\alpha[\mathcal{K}_{\infty}^{\lambda_-}(t)Qf]\right\|_{L^1}+\left\|\partial_x^\alpha[\mathcal{K}_{\infty}^{\mu_-}(t)Qf]\right\|_{L^1}\\[1ex]
 \displaystyle
 \quad
\leq 
Ct^{-(j-k)}e^{-ct}\|f\|_{W^{2k+|\alpha|,1}}.
\end{array}
\label{Sineq6} 
\end{align}
Here $a^+$ denotes $a^+=\max\{0,a\}$ for $a\in\mathbb{R}$.
\end{lem}
{\bf Proof.} We see from [20, Theorem 4.2.] that $(5.1)$--$(5.4)$ are true. Therefore it remains to show $(5.5)$ and $(5.6)$.

\vspace{1ex}
We first write $\partial_t^j\partial_x^\alpha\mathcal{K}_{\infty}^{\lambda_\pm}(t)Qf$ as
\begin{align*}
&\partial_t^j\partial_x^\alpha[\mathcal{K}_{\infty}^{\lambda_+}(t)Qf]
=\mathcal{F}^{-1}\left[\frac{ \lambda_+(\xi)^je^{\lambda_+(\xi)t}}{\lambda_+(\xi)-\lambda_-(\xi) }\hat\varphi_\infty(\xi)\frac{{\xi}^\top \xi}{|\xi|^2}\right]*\partial_x^\alpha  f.
\end{align*}

We use the formula 
\begin{equation}
e^{i\xi\cdot x}=\sum_{|\eta|=m} \frac{(-ix)^\eta}{|x|^{2m}}\partial_\xi^\eta(e^{i\xi\cdot x}). \label{eixix}
\end{equation}
By $m$-times integration by parts and $(5.7)$, we have
\begin{align*}
&\mathcal{F}^{-1}\left[\frac{ \lambda_+(\xi)^j e^{\lambda_+(\xi)t}}{\lambda_+(\xi)-\lambda_-(\xi) }\hat\varphi_\infty(\xi)\frac{{\xi}^\top \xi}{|\xi|^2}\right] \\
&=\frac{1}{(2\pi)^{\frac{3}{2}}}\sum_{|\eta|=m}\frac{(ix)^\eta}{|x|^{2m}}\int_{|\xi|\geq \sqrt{2}M_2}e^{i\xi\cdot x}\partial_\xi^\eta\left[\frac{ \lambda_+(\xi)^je^{\lambda_+(\xi)t}}{\lambda_+(\xi)-\lambda_-(\xi) }\hat\varphi_\infty(\xi)\frac{{\xi}^\top \xi}{|\xi|^2}\right]\mathrm{d}\xi.
\end{align*}
Since
\begin{align*}
&\left|\partial_\xi^\eta\left(\frac{ \lambda_+(\xi)^j e^{\lambda_+(\xi)t}}{\lambda_+(\xi)-\lambda_-(\xi) }\hat\varphi_\infty(\xi)\frac{{\xi}^\top \xi}{|\xi|^2}\right)\right|\leq C_{\eta}(1+t)^{|\eta|} e^{-c_1t}|\xi|^{-|\eta|-2},
\end{align*}
we obtain
\begin{gather}
\begin{array}{l}
\displaystyle
\left|\mathcal{F}^{-1}\left[\frac{ \lambda_+(\xi)^j e^{\lambda_+(\xi)t}}{\lambda_+(\xi)-\lambda_-(\xi) }\hat\varphi_\infty(\xi)\frac{{\xi}^\top \xi}{|\xi|^2}\right]\right| \\[2ex]
\displaystyle
\quad
\leq C_m |x|^{-m}e^{-c_1 t}\int_{|\xi|\geq\sqrt{2}M_2} |\xi|^{-m-2}\mathrm{d}\xi
\\[2ex]
\displaystyle
\quad
\leq C_m |x|^{-m}e^{-c_1 t},~m\geq2.
\end{array}
\label{estNinfeta1} 
\end{gather}
Therefore, by using $(5.8)$, we have
\begin{align*}
&\|\partial_t^j \partial_x^\alpha[\mathcal{K}_{\infty}^{\lambda_+}(t)Qf]\|_{L^1}\\
&\leq\left\|\mathcal{F}^{-1}\left[\frac{ \lambda_+(\xi)^j e^{\lambda_+(\xi)t}}{\lambda_+(\xi)-\lambda_-(\xi) }\hat\varphi_\infty(\xi)\frac{{\xi}^\top \xi}{|\xi|^2}\right]\right\|_{L^1}\|\partial_x^\alpha f\|_{L^1} \\
&\leq Ce^{-c_1t}\int_{|x|\leq 1 } |x|^{-2} \left(\int_{|\xi|\geq\sqrt{2}M_2} |\xi|^{-4} \mathrm{d}\xi\right)\mathrm{d}x\|f\|_{W^{|\alpha|,1}} \\
&\quad + Ce^{-c_1t}\int_{|x|\geq 1 } |x|^{-4} \left(\int_{|\xi|\geq\sqrt{2}M_2}|\xi|^{-6} \mathrm{d}\xi\right)\mathrm{d}x\|f\|_{W^{|\alpha|,1}}  \\
&\leq Ce^{-ct}\|f\|_{W^{|\alpha|,1}}. 
\end{align*}
Similarly, we obtain
\begin{align*}
&\|\partial_t^j \partial_x^\alpha[\mathcal{K}_{\infty}^{\mu_+}(t)Qf]\|_{L^1} \leq Ce^{-ct}\|f\|_{W^{|\alpha|,1}}.
\end{align*}
This completes the proof of $(5.5)$.

We next show $(5.6)$. 
We write $\partial_t^j\partial_x^\alpha\mathcal{K}_{\infty}^{\lambda_-}(t)Qf$ as
\begin{align*}
&\partial_t^j\partial_x^\alpha[\mathcal{K}_{\infty}^{\lambda_-}(t)Qf] \\
&=\mathcal{F}^{-1}\left[\frac{\lambda_-(\xi)^j e^{\lambda_-(\xi)t}(1+|\xi|^2)^{j-k}}{(\lambda_+(\xi)-\lambda_-(\xi))(1+|\xi|^2)^j} \hat\varphi_\infty(\xi)\frac{{\xi}^\top \xi}{|\xi|^2}\right]*\partial_x^\alpha (1-\Delta)^k f, \\
\end{align*}
where $(1-\Delta)^k=\mathcal{F}^{-1}(1+|\xi|^2)^{k}\mathcal{F}$.

By $m$-times integration by parts and $(5.7)$, we have
\begin{align*}
&\mathcal{F}^{-1}\left[\frac{ \lambda_-(\xi)^je^{\lambda_-(\xi)t}(1+|\xi|^2)^{j-k} }{(\lambda_+(\xi)-\lambda_-(\xi))(1+|\xi|^2)^j }\hat\varphi_\infty(\xi)\frac{{\xi}^\top \xi}{|\xi|^2}\right] \\
&=\frac{1}{(2\pi)^{\frac{3}{2}}}\sum_{|\eta|=m}\frac{(ix)^\eta}{|x|^{2m}}\int_{|\xi|\geq \sqrt{2}M_2}e^{i\xi\cdot x}\partial_\xi^\eta\left[\frac{ \lambda_-(\xi)^je^{\lambda_-(\xi)t}(1+|\xi|^2)^{j-k} }{(\lambda_+(\xi)-\lambda_-(\xi))(1+|\xi|^2)^j }\hat\varphi_\infty(\xi)\frac{{\xi}^\top \xi}{|\xi|^2}\right]\mathrm{d}\xi.
\end{align*}
Since
\begin{align*}
&\left|\partial_\xi^\eta\left(\frac{ \lambda_-(\xi)^je^{\lambda_-(\xi)t}(1+|\xi|^2)^{j-k}}{(\lambda_+(\xi)-\lambda_-(\xi))(1+|\xi|^2)^{j} }\hat\varphi_\infty(\xi)\frac{{\xi}^\top \xi}{|\xi|^2}\right)\right|\leq C_\eta t^{-(j-k)}e^{-c_1t-c_2|\xi|^2t}|\xi|^{-|\eta|-2},
\end{align*}
we obtain
\begin{gather}
\begin{array}{l}
\displaystyle
\left|\mathcal{F}^{-1}\left[\frac{ \lambda_-(\xi)^je^{\lambda_-(\xi)t}(1+|\xi|^2)^{j-k}}{(\lambda_+(\xi)-\lambda_-(\xi))(1+|\xi|^2)^{j} }\hat\varphi_\infty(\xi)\frac{{\xi}^\top \xi}{|\xi|^2}\right]\right| \\[1.5ex]
\displaystyle
\quad
\leq C_m|x|^{-m}t^{-(j-k)}e^{-c_1 t}\int_{|\xi|\geq\sqrt{2}M_2} e^{-c_2|\xi|^2t}|\xi|^{-m-2}\mathrm{d}\xi\\[2ex]
\displaystyle
\quad
\leq 
\begin{cases}
C_m|x|^{-m}t^{-(j-k)-\frac{3}{2}}e^{-c_1 t}, & m\geq0, \\
C_m|x|^{-m}t^{-(j-k)}e^{-c_1 t}, & m\geq2.
\end{cases}
\end{array}
\label{estNinfeta2}
\end{gather}
Therefore, by using $(5.9)$, we have
\begin{align*}
&\|\partial_t^j \partial_x^\alpha[\mathcal{K}_{\infty}^{\lambda_-}(t)Qf]\|_{L^1}\\
&\leq\left\|\mathcal{F}^{-1}\left[\frac{ \lambda_-(\xi)^j e^{\lambda_-(\xi)t}}{\lambda_+(\xi)-\lambda_-(\xi) }\hat\varphi_\infty(\xi)\frac{{\xi}^\top \xi}{|\xi|^2}\right]\right\|_{L^1}\|\partial_x^\alpha(1-\Delta)^k f\|_{L^1} \\
&\leq Ct^{-(j-k)}e^{-c_1t}\int_{|x|\leq 1 } |x|^{-2} \mathrm{d}x\|f\|_{W^{2k+|\alpha|,1}} \\
&\quad + Ct^{-(j-k)}e^{-c_1t}\int_{|x|\geq 1 } |x|^{-4} \mathrm{d}x\|f\|_{W^{2k+|\alpha|,1}}  \\
&\leq Ct^{-(j-k)}e^{-ct}\|f\|_{W^{2k+|\alpha|,1}}.
\end{align*}
Similarly, we obtain
\[
\left\|\partial_t^j\partial_x^\alpha[\mathcal{K}_{\infty}^{\mu_-}(t)Qf]\right\|_{L^1}\leq
Ct^{-(j-k)}e^{-ct}\|f\|_{W^{2k+|\alpha|,1}}.
\]
This completes the proof of $(5.6)$. $\blacksquare$

\vspace{1ex}
We next investigate the middle frequency part.
\begin{lem}\label{midfreqestetL}
The following estimates hold for $\alpha\in(\mathbb{N}\cup\{0\})^3$, $j\geq0$ and $t\geq0$:
\begin{align}
&(i)~\|\partial_t^j\partial_x^\alpha[\mathcal{K}_{M}^{\lambda_\pm}(t)f]\|_{L^1}+\left\|\partial_t^j\partial_x^\alpha[\mathcal{K}_{M}^{\mu_\pm}(t)f]\right\|_{L^1} \leq  Ce^{-ct}\|f\|_{L^1}, \label{Mineq1} 
\\[2ex] 
&(ii)~\|\partial_t^j\partial_x^\alpha[\mathcal{M}_{M}^{\lambda_\pm}(t)f]\|_{L^1}+\left\|\partial_t^j\partial_x^\alpha[\mathcal{M}_{M}^{\mu_\pm}(t)f]\right\|_{L^1} \leq  Ce^{-ct}\|f\|_{L^1}, \label{Mineq2} 
\\[2ex]
&(iii)~\|\partial_t^j\partial_x^\alpha[\mathcal{K}_{M}^{\lambda_\pm}(t)Qf]\|_{L^1}+\left\|\partial_t^j\partial_x^\alpha[\mathcal{K}_{M}^{\mu_\pm}(t)Qf]\right\|_{L^1} \leq  Ce^{-ct}\|f\|_{L^1}. \label{Mineq3} 
\end{align}
\end{lem}

{\bf Proof.} The following formulas hold for $\frac{M_1}{\sqrt{2}}\leq|\xi|\leq\sqrt{2}M_2$
\[
\frac{e^{\lambda_+(\xi)t}-e^{\lambda_-(\xi)t}}{\lambda_+(\xi)-\lambda_-(\xi)}=\frac{1}{2\pi i}\int_{\Gamma}\frac{e^{zt}}{z^2+\nu|\xi|^2z+\beta^2|\xi|^2}\mathrm{d}z,
\]
\[
\frac{e^{\mu_+(\xi)t}-e^{\mu_-(\xi)t}}{\mu_+(\xi)-\mu_-(\xi)}=\frac{1}{2\pi i}\int_{\Gamma}\frac{e^{zt}}{z^2+(\nu+\tilde{\nu})|\xi|^2z+(\beta^2+\gamma^2)|\xi|^2}\mathrm{d}z.
\]
Here $\Gamma$ is a closed path containing $\lambda_\pm(\xi)$ and $\mu_\pm(\xi)$, and included in $\{z\in\mathbb{C}|\mathrm{Re} z\leq -c_3\}$. Here $c_3$ is a positive number taken by
\[
\max_{\frac{M_1}{\sqrt{2}}\leq|\xi|\leq\sqrt{2}M_2}\mathrm{Re}\mu_j(\xi)\leq-2c_3, j=1,2,3,4. 
\]
Hence, we can compute as in a similar manner to [11,20] to obtain
\begin{align*}
&
\left|\partial_t\partial_x^\alpha\mathcal{F}^{-1}\left[\frac{1}{2\pi i}\int_{\Gamma}\frac{e^{zt}}{z^2+\nu|\xi|^2z+\beta^2|\xi|^2}\mathrm{d}z\eta(\xi)\hat\varphi_M(\xi)\right]\right|\\
&\leq C_{j,\alpha,N}e^{-ct}|x|^{-N}, j+|\alpha|\geq1, N\geq0,
\end{align*}
\begin{align*}
&\left|\partial_t\partial_x^\alpha\mathcal{F}^{-1}\left[\frac{1}{2\pi i}\int_{\Gamma}\frac{e^{zt}}{z^2+(\nu+\tilde{\nu})|\xi|^2z+(\beta^2+\gamma^2)|\xi|^2}\mathrm{d}z\eta(\xi)\hat\varphi_M(\xi)\right]\right|\\
&\leq C_{j,\alpha,N}e^{-ct}|x|^{-N}, j+|\alpha|\geq1, N\geq0.
\end{align*}
Here $\eta$ is a function such that $\eta\in C^\infty(S^2),~S^2=\{\xi\in\mathbb{R}^3~|~|\xi|=1\}$ and $\eta=\eta(\frac{\xi}{|\xi|})$.
Therefore we have $(5.10)$--$(5.12)$.
This completes the proof. $\blacksquare$ \\[2ex]
{\bf Proof of Proposition \ref{L1estetL} (ii).} 
By using Lemma 5.1 and Lemma 5.2, we have the following estimates of $\mathcal{K}^{j1}(t)P_\infty\phi_0,~j=1,2,$ and $\mathcal{K}^{j2}(t)P_\infty w_0,~j=1,2,3$
\begin{align*}
\|\mathcal{K}^{11}(t)P_\infty\phi_0\|_{L^1}&\leq Ce^{-ct}\|\phi_0\|_{L^1},\\[0.5ex]
\|\mathcal{K}^{12}(t)P_\infty w_0\|_{L^1}&\leq Ce^{-ct}\|w_0\|_{L^1},\\[0.5ex]
\|\mathcal{K}^{21}(t)P_\infty\phi_0\|_{L^1}&\leq Ce^{-ct}\|\phi_0\|_{L^1}, \\[0.5ex]
\|\mathcal{K}^{22}(t)P_\infty w_0\|_{L^1}&\leq 
Ce^{-ct}\|w_0\|_{W^{2,1}},
 \\[0.5ex]
\|\mathcal{K}^{32}(t)P_\infty w_0\|_{L^1}&\leq 
Ce^{-ct}\|w_0\|_{W^{2,1}},
\end{align*}
We next focus on $\mathcal{K}^{j3}(t)P_\infty G_0,~j=2,3.$ By using $G_0=\nabla\psi_0$, we write $\mathcal{F}[\mathcal{K}^{j3}(t) G_0],~j=2,3$ as
\begin{align*}
&\mathcal{F}[\mathcal{K}^{23}(t)G_0](\xi) \\[0.5ex]
&=-i\beta^2\frac{e^{\lambda_+(\xi)t}-e^{\lambda_-(\xi)t}}{\lambda_+(\xi)-\lambda_-(\xi)}(\hat{ G}_0(\xi)-\mathrm{tr}(\hat{G}_0(\xi))I)\xi\  \\[0.5ex]
&\quad\quad-i\beta^2\frac{e^{\mu_+(\xi)t}-e^{\mu_-(\xi)t}}{\mu_+(\xi)-\mu_-(\xi)}\mathrm{tr}(\hat{G}_0(\xi))\xi,\\[0.5ex]
&\mathcal{F}[\mathcal{K}^{33}(t)G_0](x) \\[0.5ex]
&=\frac{\lambda_+(\xi)e^{\lambda_+(\xi)t}-\lambda_-(\xi)e^{\lambda_-(\xi)t}}{\lambda_+(\xi)-\lambda_-(\xi)}\hat{G}_0(\xi)+\nu|\xi|^2\frac{e^{\lambda_-(\xi)t}-e^{\lambda_+(\xi)t}}{\lambda_+(\xi)-\lambda_-(\xi)}\hat{G}_0(\xi)  \\[0.5ex]
&+\left((\nu+\tilde{\nu})\frac{e^{\mu_+(\xi)t}-e^{\mu_-(\xi)t}}{\mu_+(\xi)-\mu_-(\xi)}-\nu\frac{e^{\lambda_-(\xi)t}-e^{\lambda_+(\xi)t}}{\lambda_+(\xi)-\lambda_-(\xi)}\right)\xi{}^\top\xi G_0(\xi)  \\[0.5ex]
&\quad\quad+\left(\frac{\mu_+(\xi)e^{\mu_+(\xi)t}-\mu_-(\xi)e^{\mu_-(\xi)t}}{\mu_+(\xi)-\mu_-(\xi)}-\frac{\lambda_+(\xi)e^{\lambda_+(\xi)t}-\lambda_-(\xi)e^{\lambda_-(\xi)t}}{\lambda_+(\xi)-\lambda_-(\xi)}\right)\frac{\xi^{\top}\xi}{|\xi|^2}\hat{G}_0(\xi).
\end{align*}
It then follows from Lemma 5.1 and Lemma 5.2 that
\begin{align*}
\|\mathcal{K}^{23}(t)P_\infty G_0\|_{L^1}&\leq 
Ce^{-ct}\|G_0\|_{W^{2,1}},\\
\|\mathcal{K}^{33}(t)P_\infty G_0\|_{L^1}&\leq 
Ce^{-ct}\|G_0\|_{W^{2,1}}.
\end{align*}
Consequently we arrive at
\[
\|e^{-tL}P_\infty U_0\|_{L^1}\leq 
Ce^{-ct}\|U_0\|_{W^{2,1}}.
\]
This completes the proof. $\blacksquare$

\noindent
\textbf{Acknowledgements. } 
This work was partially supported by JSPS KAKENHI Grant Number 19J10056. 

\bibliographystyle{unsrt}  

\end{document}